\documentclass{amsart}
\usepackage{graphicx}
\usepackage{amssymb}
\newtheorem{theorem}{Theorem}[section]

\newtheorem{corollary}[theorem]{Corollary}

\theoremstyle{definition}

\newtheorem{example}[theorem]{Example}
\newtheorem{remark}[theorem]{Remark}

\newcommand{\R}{{\bf{R}}}
\newcommand{\C}{{\bf{C}}}
\newcommand{\Z}{{\bf{Z}}}
\newcommand{\Q}{{\bf{Q}}}
\newcommand{\D}{\Delta}
\newcommand{\la}{\langle}
\newcommand{\ra}{\rangle}

\newcommand{\Sg}{\Sigma_g}
\newcommand{\Mg}{\mathcal{M}_g}
\newcommand{\Ig}{\mathcal{I}_g}
\newcommand{\Kg}{\mathcal{K}_g}
\newcommand{\Hg}{\mathcal{H}_g}

\newcommand{\G}{\Gamma}

\begin{document}

\thanks{Primary 57R20; Secondary 57M27.}

\thanks{To appear in the Proceedings of the Postnikov Memorial Conference.}

\title[On a secondary invariant]{On a secondary invariant of the \\
hyperelliptic mapping class group}

\author{Takayuki Morifuji}

\address{Department of Mathematics, Tokyo University of 
Agriculture and Technology, 
2-24-16 Naka-cho, Koganei, Tokyo 184-8588, Japan}

\email{morifuji@cc.tuat.ac.jp}

\maketitle

\begin{abstract}In this paper, 
we discuss relations among several invariants 
of $3$-manifolds including Meyer's function, 
the $\eta$-invariant, the von Neumann $\rho$-invariant 
and the Casson invariant 
from the viewpoint of the  mapping class group of a surface. 
\end{abstract}

\section{Introduction}
 
Let $\Sg$ be a closed oriented smooth surface of genus $g$ 
and $\mathrm{Diff}_+\Sg$ the group of 
orientation preserving diffeomorphisms of $\Sg$ 
equipped with the $C^\infty$-topology. 
The mapping class group $\Mg$ is defined to be 
the group of path components of 
$\mathrm{Diff}_+\Sg$. 
In this paper, 
by an automorphism of $\Sg$, 
we mean an element of $\Mg$. 
Next 
let 
$r:\Mg\to\mathrm{Sp}(2g,\Z)$ 
denote the classical representation defined 
by the action of $\Mg$ on the first integral homology group 
$H_1(\Sg,\Z)$. 
It is known that 
$r$ is surjective and 
the Torelli group $\Ig$ is defined to be 
the kernel of $r$. 
Namely, 
it has the trivial action on $H_1(\Sg,\Z)$. 

Our main object here is 
the signature cocycle $\tau$, 
which is a group $2$-cocycle of the Siegel modular group 
$\mathrm{Sp}(2g,\Z)$. 
It was introduced by Meyer (see \cite{[Me]}) 
to describe a signature formula for surface bundles 
over a surface. 
If 
we pull back the cocycle by the representation $r$, 
then 
we can think of $\tau$ as a $2$-cocycle of $\Mg$. 
But here, 
we mainly consider the restriction of $\tau$ 
to a subgroup of $\Mg$, 
that is, 
the hyperelliptic mapping class group 
$\D_g$, 
which is the centralizer with 
a fixed hyperelliptic involution. 
As is known, 
$\D_g=\Mg$ if $g=1,2$ and 
$\D_g\not=\Mg$ if $g\geq 3$. 
Moreover 
the fact that 
$H^*(\D_g,\Q)=0$ for $*=1,2$ 
(see \cite{[C]}, \cite{[Kw]}) 
implies 
the restriction of $\tau$ to $\D_g$ 
must be the coboundary of the unique rational 
$1$-cochain $\phi:\D_g\to\Q$ 
(i.e. $\delta\phi=r^*\tau|_{\D_g}$ holds). 
In the following, 
we call it Meyer's function 
of the hyperelliptic mapping class group. 
Namely, 
it is a secondary invariant of $\D_g$ 
associated with the signature cocycle. 

In the case of genus one, 
explicit formulas of $\tau$ and $\phi$ 
were given by Meyer himself, 
Kirby-Melvin \cite{[KM]} and Sczech \cite{[Sc]}. 
On the other hand, 
geometric aspects of Meyer's function 
have been studied by Atiyah in \cite{[A1]}. 
In fact, 
he related it to various invariants 
defined for each element of 
$\D_1=\mathcal{M}_1\cong SL(2,\Z)$ 
including 
Hirzebruch's signature defect, 
the logarithmic monodromy of Quillen's 
determinant line bundle, 
the Atiyah-Patodi-Singer $\eta$-invariant 
and its adiabatic limit, 
the special value of the Shimizu $L$-function 
and so on. 

The purpose of the present paper 
is to discuss some extensions of Atiyah's results 
to higher genus case 
and 
explain relations among several 
invariants of $3$-manifolds, 
that is, 
Meyer's function, 
the $\eta$-invariant, 
the von Neumann $\rho$-invariant 
and the Casson invariant 
from the viewpoint of the mapping class group 
of a surface. 

Now 
we describe the contents of this paper. 
In the next section, 
we recall the definitions of 
the signature cocycle and Meyer's function 
according to \cite{[Me]} and \cite{[Mf5]}. 
In Section 3 
we give a relation between Meyer's function 
and the $\eta$-invariant under 
the mapping torus construction. 
In Section 4 
we describe the von Neumann $\rho$-invariant 
for a $\Z$-covering by using Meyer's function. 
In the last section, 
the Casson invariant of an integral homology 
$3$-sphere is described by Meyer's function 
through the correspondence between 
elements of $\Mg$ and $3$-manifolds 
via Heegaard splittings. 

\section{Signature cocycle and Meyer's function}
 
In this section, 
we quickly review the definitions of the signature cocycle and 
Meyer's function. 
First 
we take two elements 
$A,B\in \mathrm{Sp}(2g,\Z)$ 
and 
define the vector space $V_{A,B}$ by 
$$
V_{A,B}
=
\left\{
(x,y)\in \R^{2g}\times\R^{2g}
~|~(A^{-1}-I)x+(B-I)y=0
\right\},
$$
where 
$I$ denotes the identity matrix. 
We also define the pairing map on $\R^{2g}\times\R^{2g}$ by 
$$
\la
(x_1,y_1),(x_2,y_2) 
\ra_{A,B}
=(x_1+y_1)\cdot J(I-B)y_2,
$$
where 
$\cdot$ denotes the standard 
inner product in $\R^{2g}$ 
and 
$J$ is the $2g\times 2g$ real matrix 
corresponding to the multiplication 
by $\sqrt{-1}$ on $\C^g=\R^{2g}$. 
We can easily check that 
this pairing map is 
a symmetric bilinear form on 
$V_{A,B}$ (possibly degenerate), 
and so 
we define $\tau(A,B)$ to be 
the signature of 
the pairing map 
$\la~,~\ra_{A,B}$ on 
the vector space $V_{A,B}$. 

From Novikov additivity, 
$\tau(A,B)$ satisfies 
so called the cocycle condition, 
namely 
$$
\tau(A,B)+\tau(AB,C)
=
\tau(A,BC)+\tau(B,C),
$$
so that 
$\tau$ is a group $2$-cocycle 
of the Siegel modular group 
and represents the cohomology class 
$-4c_1\in H^2(\mathrm{Sp}(2g,\Z),\Z)$. 
We call this $2$-cocycle 
$\tau$ the signature cocycle \cite{[Me]}. 

\begin{remark}
(i) 
The signature cocycle $\tau$ is 
defined for $\mathrm{Sp}(2g,\Z)$, 
but 
we can regard it as a $2$-cocycle of $\Mg$ in terms of 
the homology representation $r$. 
We use the same letter $\tau$ for simplicity. 
(ii) 
$\tau(A,B)$ just presents $\sigma(W)$, 
the signature of the $4$-manifold $W$, 
which is a surface bundle over the pair of pants 
with three boundary components 
$M_A,~M_B$ and $-M_{AB}$, 
where 
$M_A$ denotes a mapping torus with 
monodromy $A$. 
(iii) 
By definition, 
$\tau$ is a bounded $2$-cocycle. 
In fact, 
the range of $|\tau|$ is clearly 
bounded by $2g$. 
\end{remark}

Now 
in the following, 
we consider the hyperelliptic mapping class group 
of a surface $\Sg$. 
We fix an involution $\iota\in\Mg$ 
with $2g+2$ fixed points. 
The centralizer 
$$
\D_g
=\{f\in\Mg~|~f\iota=\iota f\}
$$
is called the hyperelliptic mapping class group 
of $\Sg$ with respect to $\iota$. 
The conjugacy class of $\D_g$ 
in $\Mg$ does not depend on a choice of $\iota$. 

As was mentioned in the introduction, 
the rational cohomology groups of $\D_g$ 
are trivial for dimensions $1$ and $2$. 
Hence 
the cohomology class represented by $\tau$ 
is a torsion element in $H^2(\D_g,\Z)$. 
More precisely, 
we can show that the order of $[\tau]$ in 
$H^2(\D_g,\Z)$ is $2g+1$ (see \cite{[E]}). 
Thus 
there exists the uniquely defined mapping 
$$
\phi:\D_g\to\frac{1}{2g+1}\Z
$$ 
so that 
the coboundary of $\phi$ coincides with 
the restriction of $\tau$ to $\D_g$, 
where 
$({1}/({2g+1}))\Z$ denotes 
an additive group 
$\{{n}/({2g+1})\in\Q~|~n\in\Z\}$. 
Namely, 
we have 
$$
\delta\phi(f_1,f_2)
=\phi(f_2)-\phi(f_1f_2)+\phi(f_1)
=\tau(f_1,f_2)
$$ 
for $f_1,f_2\in\D_g$. 
We call it Meyer's function of 
the hyperelliptic mapping class group. 

\begin{remark}
(i) Meyer's function $\phi$ is a class function 
of $\D_g$. 
That is, 
for two elements $f_1,f_2\in\D_g$, 
$\phi(f_2f_1f_2^{-1})=\phi(f_1)$ holds. 
Consequently, 
we can regard $\phi$ as an invariant 
of surface bundles over the circle. 
(ii) The relation 
$\delta\phi=\tau|_{\D_g}~(g\geq2)$ 
implies that 
$\phi$ is a homomorphism on 
the Torelli group $\Ig\cap\D_g$, 
because 
$\tau$ is originally defined 
on the Siegel modular group 
$\mathrm{Sp}(2g,\Z)$. 
\end{remark}

The group presentation of $\D_g$ was 
given by Birman and Hilden (see \cite{[BH]}). 
By using its defining relations, 
we can determine the values of $\phi$ 
on generators of $\D_g$. 

\begin{example} 
For the generators 
$\zeta_1,\zeta_2,\ldots,\zeta_{2g+1}$ 
of $\D_g~(g\geq 2)$, 
the values of Meyer's function $\phi$ 
are equal to ${(g+1)}/{(2g+1)}$. 
This is shown as follows. 
First 
each generator is mutually conjugate 
(in fact, 
$\zeta_{i+1}=\xi\zeta_i\xi^{-1}$ 
holds for 
$\xi=\zeta_1\cdots\zeta_{2g+1}$), 
so that 
$\phi$ has the same value on generators. 
Moreover 
substituting a defining relation 
of $\D_g$ for $\delta\phi=\tau$, 
we can evaluate $\phi(\zeta_i)$ 
explicitly. 
\end{example}

\begin{example}
Let 
$\psi_h\in\D_g~(1\leq h\leq g-1)$ 
be a BSCC-map of genus $h$. 
Namely, 
it is a Dehn twist along a bounding simple closed curve 
on $\Sg$ which is invariant under the action of 
the hyperelliptic involution $\iota$ and 
separates $\Sg$ into two subsurfaces 
of genus $h$ and $g-h$. 
Then 
$\psi_h$ is presented by 
$\psi_h=f(\zeta_1\cdots\zeta_{2h})^{4h+2}f^{-1}$ 
for some automorphism $f\in\D_g$. 
Thus we have 
$\phi(\psi_h)=-4h(g-h)/(2g+1)$. 
\end{example}

\begin{remark}
(i) There exist other constructions of Meyer's function 
of genus $2$, which are due to Kasagawa \cite{[Kg]} and 
Iida \cite{[I]}. 
(ii) Meyer's function plays an important role 
in the study of certain singular fibrations. 
In fact, 
Matsumoto \cite{[Ma]} and Endo \cite{[E]} defined the local 
signature of hyperelliptic Lefschetz 
fibrations by using 
Meyer's function. 
See also \cite{[Ku]} for Meyer's function of 
plane curves. 
\end{remark}

\section{Eta-invariant}
 
In this section, 
we discuss a relation between Meyer's function 
and the $\eta$-invariant of 
a surface bundle over the circle. 
Let 
$M$ be an oriented closed Riemannian $3$-manifold. 
By using the spectrum of the signature operator, 
the $\eta$-invariant of $M$ is defined 
to be the value at zero of the $\eta$-function 
associated with the signature operator. 
Instead of giving the precise definition of 
the $\eta$-invariant, 
we recall the following index theorem 
due to Atiyah, Patodi and Singer 
(see \cite{[APS]}). 
Let 
$W$ be a compact oriented Riemannian $4$-manifold 
with product metric near the boundary 
$\partial W=M$. 
Then 
the $\eta$-invariant of $M$ is given by 
$$
\eta(M)
=\frac13\int_Wp_1-\sigma(W),
$$
where 
$p_1$ is the first Pontrjagin form 
for the Riemannian connection on 
the tangent bundle $TW$ of $W$. 
As is well-known, 
if $W$ is a closed $4$-manifold, 
then 
the signature $\sigma(W)$ is described by 
the integral of $p_1$. 
That is, 
the above formula gives 
the extension of the Hirzebruch signature formula 
to $4$-manifold with boundary. 
Here 
we have to remark that 
the $\eta$-invariant is not a topological 
invariant, but a spectral invariant. 

Now 
for an automorphism $f\in \Mg$, 
let $M_f$ be the mapping torus 
corresponding to $f$. 
Namely, 
it is the identification space 
$\Sg\times\R/(x,t)\sim(f(x),t+1)$. 
When 
$f\in\Mg$ 
is of finite order (or a periodic automorphism), 
we endow $M_f$ with the metric 
which is induced from the product 
of the standard metric on the circle and 
an $f$-invariant metric of $\Sg$. 
If 
we restrict ourselves to 
the hyperelliptic mapping class group $\D_g$, 
we obtain the following theorem. 

\begin{theorem}
$\eta(M_f)=\phi(f)$ holds for any periodic automorphism 
$f\in\D_g$.
\end{theorem}

This theorem 
is shown by using an explicit formula 
of the $\eta$-invariant for a periodic automorphism 
of the mapping class group $\Mg$ (see \cite{[Mf1]}). 
More precisely, 
it holds that 
$$
\eta(M_f)
=
\frac1n\sum_{k=1}^{n-1}\tau\left(f,f^k\right)
$$
for $f\in\Mg$ of the order $n$. 
The formula is based on the index theorem, 
due to Atiyah, Patodi and Singer. 
The point is to use the notion of 
$2$-framing of $3$-manifolds, 
which 
is a trivialization of twice the tangent bundle 
$2TM=TM\oplus TM$ 
as a $\mathrm{Spin}(6)$-bundle 
(see \cite{[A2]} for details). 

\begin{example} 
We can explicitly calculate all the values of 
Meyer's function $\phi$ 
on periodic automorphisms for the genus $1$ and $2$ 
cases by using the homological action of each 
periodic automorphism 
(see \cite{[Mf2]}). 
In particular, 
we notice a relation to 
the Nielsen-Thurston theory of surface automorphisms \cite{[CB]}. 
As is known, 
the theory classifies the automorphisms of $\Sg$ 
into the following three types: 
(i) periodic, 
(ii) reducible and 
(iii) pseudo-Anosov. 
We 
say an automorphism of $\Sg$ is reducible, 
if it has a representative which leaves some essential 
$1$-submanifold of the surface invariant. 
We then easily see that 
(i) and (ii) have some overlap, 
although (iii) cannot have any intersection 
with others. 
Hence, 
a given periodic automorphism 
is either reducible or irreducible. 
In fact, 
there are several works concerning 
the characterization of the reducibility 
of periodic automorphisms 
(see \cite{[Kh]}, \cite{[Mf3]} and their references). 

From direct computations of 
Meyer's function, 
if the genus of $\Sg$ is $1$ or $2$, 
then 
we obtain the following. 
That is, 
a periodic automorphism $f$ is reducible if and only if 
$\phi(f)=0$ (i.e. $\eta(M_f)=0$). 
Thus we can say that 
the reducibility of periodic automorphisms of 
a surface with low genus is characterized by 
the vanishing of Meyer's function 
or the $\eta$-invariant. 
In general, 
the same statement as above does not hold 
for higher genus case, 
but 
we can get the similar result for automorphisms 
of prime orders of 
hyperelliptic Riemann surfaces (see \cite{[Mf4]}). 
\end{example}

Now 
as an immediate corollary of Theorem 3.1, 
we have the next assertion. 

\begin{corollary}
Let 
$f\in\Mg$ be a periodic automorphism. 
If $f\in\D_g$, 
that is, 
$f$ commutes with the fixed hyperelliptic involution $\iota$, 
then $\eta(M_f)\in ({1}/({2g+1}))\Z$ 
holds. 
\end{corollary}

Thus 
in some sense, 
we may regard the $\eta$-invariant 
as an obstruction for distinguishing 
a given periodic automorphism to be 
hyperelliptic or not. 

\begin{example}
Let 
$f\in\mathcal{M}_3$ be a periodic automorphism of order $3$ 
so that 
the quotient orbifold of $\Sigma_3$ 
by its cyclic action is homeomorphic to 
$S^2(3,3,3,3,3)$ 
($2$-sphere with five cone points of index $3$). 
Then an easy computation shows that 
the $\eta$-invariant of corresponding mapping torus 
is given by 
$$
\eta(M_f)=-\frac23\not\in \frac17\Z.
$$
Hence, 
Corollary 3.3 implies that 
$f$ cannot commute with the hyperelliptic 
involution $\iota$. 
\end{example}

\begin{remark}
{Let $G$ be the finite cyclic group generated by 
a periodic automorphism $f\in\Mg$. 
Since 
$H^*(G,\Q)=0$ for $*=1,2$, 
we can define Meyer's function 
$\phi_G:G\to\Q$ of $G$. 
Thereby 
as an analogue of Theorem 3.1, 
we see that 
$\eta(M_f)=\phi_G(f)$ holds 
for any periodic automorphism 
$f\in\Mg$ (see \cite{[Ak]}). 
}
\end{remark}

\section{von Neumann rho-invariant} 

Next 
we discuss a relation to the von Neumann $\rho$-invariant. 
Let 
$\G$ be a discrete group and 
$M$ an oriented closed Riemannian $3$-manifold. 
Moreover, 
we assume that 
we are given a surjective homomorphism 
from the fundamental group 
$\pi_1M$ to $\G$. 
Then 
we can take a $\G$-covering 
$\G\to \hat{M}\to M$. 
Lifting 
the metric and the signature operator 
to $\hat{M}$, 
the $\eta$-invariant of $\hat{M}$, 
denoted by $\eta^{(2)}(\hat{M})$, 
is defined. 
We call it 
the von Neumann $\eta$-invariant. 
It is also known that 
$\eta^{(2)}(\hat{M})$ fits an index theorem as before. 

The $\eta$-invariant and the von Neumann 
$\eta$-invariant depend on the Riemannian metric, 
but 
Cheeger-Gromov showed in \cite{[CG]} 
that 
the difference of them is independent of 
the Riemannian metric. 
We denote the difference 
$\eta^{(2)}(\hat{M})-\eta(M)$ 
by 
$\rho^{(2)}(\hat{M})$ and 
call it the von Neumann $\rho$-invariant. 

In a certain sense, 
this invariant $\rho^{(2)}$ 
is an extension of the classical 
$\rho$-invariant 
due to Atiyah-Patodi-Singer. 
That is, 
the difference between the $\eta$-invariant 
$\eta_\gamma$ twisted by 
a unitary representation 
$\gamma:\pi_1M\to U(n)$ 
and the original $\eta$-invariant 
does not depend on the choice of 
the Riemannian metric. 

For 
an automorphism $f\in\D_g$, 
we consider the $\Z$-covering 
$\Z\to \hat{M}_f\to M_f$ 
associated with a surjective homomorphism 
$\pi_1M_f\to\pi_1S^1\cong\Z$. 
A simple observation tells us 
that Meyer's function is not multiplicative 
for coverings. 
But if we take a limit, 
we see that it is related to the von Neumann $\rho$-invariant. 
Namely, 
we have 

\begin{theorem}
$\displaystyle{
\rho^{(2)}(\hat{M}_f)
=\lim_{k\to\infty}
\frac{\phi(f^k)-k\phi(f)}{k}.
}$
\end{theorem}

To be more precise, 
the von Neumann $\rho$-invariant 
$\rho^{(2)}(\hat{M}_f)$ 
is given by the limit 
of the deviation from the multiplicativity 
of $\phi$ for a finite covering 
$M_{f^k}\to M_f$. 
This is shown by using 
Propositions 2.2 and 3.1 in \cite{[Mf1]} 
and the approximation theorem of 
the $\eta$-invariant, 
due to 
Vaillant \cite{[V]} and L\"{u}ck-Schick \cite{[LS]}, 
which states that 
$$
\eta^{(2)}(\hat{M})
=
\lim_{k\to\infty}
\frac{\eta(M_{(k)})}{[\G:\G_k]}
$$ 
holds 
for a descending sequence of normal subgroups 
$\G
\supset
\G_1
\supset
\G_2\supset
\cdots$ 
such that 
$[\G:\G_k]<\infty$ and $\cap_k\G_k=\{1\}$, 
and 
$\G/\G_k$-covering 
$M_{(k)}=\hat{M}/\G_k\to M$.

\begin{example} 
As an example, 
we consider the genus one case. 
In this case, 
an element $A\in SL(2,\Z)$ is 
classified into the following three cases: 

(i) Elliptic case (namely, $|\mathrm{tr}~A|<2$). 
Let $A_n\in SL(2,\Z)$ have the order $n~(n=3,4,6)$. 
We can take 
$$
A_3=
\left(
\begin{array}{cc}
-1&-1\\
1&0
\end{array}
\right),\quad
A_4=
\left(
\begin{array}{cc}
0&-1\\
1&0
\end{array}
\right)\quad
\mathrm{and}
\quad
A_6=
\left(
\begin{array}{cc}
0&-1\\
1&1
\end{array}
\right).
$$ 
An easy calculation shows that 
$$
\phi({A_3})=-\phi({A_3}^2)=-2/3,~
\phi({A_4})=-\phi({A_4}^3)=-1,~
\phi({A_4}^2)=0,
$$ 
$$
\phi({A_6})=-\phi({A_6}^5)=-4/3,~
\phi({A_6}^2)=-\phi({A_6}^4)=-2/3,~
\phi({A_6}^3)=0.
$$ 
Needless to say, 
$\phi({A_n}^n)=\phi(I)=0$ holds. 
Hence 
we have 
$$
\rho^{(2)}(\hat{M}_{A_n})
=
\left\{
\begin{array}{ll}
2/3,&n=3,\\
1,&n=4,\\
4/3,&n=6
\end{array}
\right.
$$
in terms of Theorem 4.1. 
It should be noted that 
$\rho^{(2)}(\hat{M}_f)=0$ 
for any involution $f\in\Mg$ 
(see \cite{[Mf1]}, \cite{[Mf6]}). 

(ii) Parabolic case (namely, $|\mathrm{tr}~A|=2$). 
We can take 
$\displaystyle{
A_b=
\left(
\begin{array}{cc}
1&b\\
0&1
\end{array}
\right)~
(b\in\Z).
}$
Then 
we obtain 
$\rho^{(2)}(\hat{M}_{A_b})=-\mathrm{sgn}(b)$, 
where 
$\mathrm{sgn}(b)=b/|b|$ if $b\not=0$ 
and $0$ if $b=0$. 
This result follows from Atiyah's calculation 
in \cite{[A1]}. 

(iii) 
Hyperbolic case 
(namely, $|\mathrm{tr}~A|>2$). 
Since Meyer's function of genus one satisfies 
$\phi(A^k)=k\phi(A)$ 
for a hyperbolic element $A\in SL(2,\Z)$ 
(see \cite{[Me]}), we have 
$\rho^{(2)}(\hat{M}_A)=0$ 
by virtue of Theorem 4.1. 
\end{example}

As an immediate corollary of Theorem 4.1 
and Remark 2.2 (ii), 
we obtain the following for higher genus case. 

\begin{corollary}
If $f$ is an automorphism in $\Ig\cap\D_g$, 
then 
$\rho^{(2)}(\hat{M}_f)=0$. 
\end{corollary}

Here 
we have a remark. 
If we restrict the above theorem 
to the level $2$ subgroup of $\Mg$, 
we can describe a relation 
among the von Neumann $\rho$-invariant, 
the first Morita-Mumford class $e_1\in H^2(\Mg,\Z)$ 
(see \cite{[Mt1]}, \cite{[Mu]} and Section 5 below) 
and the Rochlin invariant of a spin $3$-manifold 
in a framework 
of the bounded cohomology \cite{[G]}. 
Roughly speaking, 
the pull-back of $e_1$ into 
$H^2_b(S^1,\Z)\cong\R/\Z$ 
via a holonomy homomorphism 
is given by 
the linear combination of 
the Rochlin invariant and 
the von Neumann $\rho$-invariant 
(see \cite{[Ki2]}, \cite{[Mf6]} for details). 
We 
describe it more explicitly in the appendix 
of the present paper. 

\section{Casson invariant} 

In this last section, 
we explain a relation between Meyer's function 
and the Casson invariant. 
The Casson invariant $\lambda(M)$ 
is an integer valued invariant 
defined for an oriented integral homology $3$-sphere $M$. 
Roughly speaking, 
it counts the number of conjugacy classes 
of irreducible representations of the 
fundamental group $\pi_1M$ 
into the Lie group $SU(2)$. 

On the other hand, 
from the theory of characteristic classes of surface bundles, 
due to Morita (see \cite{[Mt2]}, \cite{[Mt3]}), 
the Casson invariant $\lambda$ can be interpreted as 
a secondary invariant associated with 
the first Morita-Mumford class $e_1$, 
through the correspondence between 
elements of the mapping class group 
and $3$-manifolds via Heegaard splittings. 

Let 
us state more precisely. 
We define $\Kg$ to be the subgroup of $\Mg$ 
generated by Dehn twists along bounding simple closed curves 
on $\Sg$. 
This is also a subgroup of the Torelli group $\Ig$, 
because clearly, 
the action of a generator of $\Kg$ 
on the homology of $\Sg$ is trivial. 
Next 
we fix a Heegaard splitting of the $3$-sphere 
$S^3=\Hg\cup_{\iota_g}-\Hg$, 
where 
$\Hg$ denotes the handle body of genus $g$ and 
$\iota_g$ is the gluing map of this splitting. 

For an automorphism $f\in \Kg$, 
we construct a $3$-manifold $M^f$ 
by regluing two handle bodies via 
the composition $\iota_g f$. 
Then 
it is easy to see that 
the resulting manifold $M^f$ is again 
a homology $3$-sphere, 
because 
$f$ acts on 
$H_1(\Sg,\Z)$ trivially. 
Thereby 
we can evaluate its Casson invariant 
$\lambda(M^f)$. 
In short, 
Morita's results claim that 
there exists a homomorphism 
$\lambda^*:\Kg\to\Z$ such that 
$\lambda^*(f)=\lambda(M^f)$. 
More precisely, 
$\lambda^*$ consists of 
the sum of two homomorphisms. 
One is Morita's homomorphism 
$d_0:\Kg\to\Q$, 
which is the core of the Casson invariant 
from the viewpoint of the mapping class group. 
The other is 
the Johnson homomorphism, 
which is interpreted as Massey higher products of 
mapping tori (see \cite{[Ki1]}). 

Now 
assume $g\geq 2$. 
We then have 

\begin{theorem}
Meyer's function essentially coincides with Morita's 
homomorphism on $\Kg\cap\D_g$. 
To be more precise, 
$$
\phi(f)=\frac13 d_0(f)
$$
holds for any automorphism $f\in\Kg\cap\D_g$. 
\end{theorem}

Therefore, 
in principle, 
we can say that the Casson invariant of an integral 
homology $3$-sphere is determined by Meyer's 
function. 
Here, 
as an example, 
we evaluate Morita's homomorphism on 
a typical element in $\Kg\cap\D_g$. 

\begin{example} 
Let $\psi_h\in\Kg\cap\D_g$ be a BSCC-map of genus $h$. 
Then 
the value of Morita's homomorphism on $\psi_h$ 
is given by 
$$
d_0(\psi_h)
=-\frac{12}{2g+1}h(g-h)
$$
in terms of Theorem 5.1 and Example 2.4. 
\end{example}

Finally 
let us review the definition of Morita's homomorphism 
very briefly. 
It is a secondary invariant associated with 
the first Morita-Mumford class $e_1\in H^2(\Mg,\Z)$. 

Let 
$e\in H^2(\Mg{}_{,*},\Z)$ be the Euler class of 
the central $\Z$-extension 
$\Z\to\Mg{}_{,1}\to\Mg{}_{,*}
$,
where 
$\Mg{}_{,1}$ and $\Mg{}_{,*}$ 
denote the mapping class groups of $\Sg$ relative to 
an embedded disc $D\subset \Sg$ and 
a base point $*\in D$ respectively. 
The center $\Z$ is generated by 
the Dehn twist parallel to $\partial D$. 
We define $e_1$ to be the Gysin image 
(integration along the fiber) of $e^2$. 
This is an element of $H^2(\Mg,\Z)$ and 
called the first Morita-Mumford class. 

Now 
it is known that there exist two canonical 
$2$-cocycles over $\Q$ representing $e_1$. 
One is the signature cocycle $\tau$ 
($e_1=[-3\tau]$ holds) 
and 
the other is the intersection cocycle $c$ 
(the latter one is defined 
once we fix a certain crossed 
homomorphism of $\Mg$). 
Therefore, 
we have the uniquely defined mapping 
$d:\Mg\to\Q$ so that 
$\delta d=c+3\tau$. 
The uniqueness follows from the fact that 
$\Mg$ is perfect for $g\geq3$. 
We denote the restriction of $d$ to 
the subgroup $\Kg$ by $d_0$. 
Then 
Morita showed that 
$d_0$ does not depend on the choice of 
the crossed homomorphism 
and 
it serves a generator of 
$H^1(\Kg,\Z)^{\Mg}$. 
This $d_0$ is the map, 
which we have called 
Morita's homomorphism. 

The point of the proof of Theorem 5.1 
is to construct a crossed homomorphism of 
$\Mg$ so that the restriction of it to 
the hyperelliptic mapping class group 
$\D_g$ is zero map. 
Using it 
to define the intersection cocycle, 
we can obtain a relation between 
Meyer's function and Morita's 
homomorphism. 

\begin{remark}
{In \cite{[Mt4]}, 
Morita gave an interpretation of $d_0$ 
in terms of Hirzebruch's signature defect 
of certain framed $3$-manifolds. 
Combining it and Theorem 5.1, 
we have a generalization of 
Atiyah's result mentioned in the introduction. 
}
\end{remark}

\section{Appendix}
 
We give here a description of the first Morita-Mumford class $e_1$ 
via the Rochlin invariant and the von Neumann $\rho$-invariant 
in the bounded cohomology $H_b^*$. 

First, 
we note that 
$e_1$ is a bounded cohomology class (see Remark 2.1 (iii)) and 
consider $e_1$ on $\Mg{}_{,*}$ rather than on $\Mg$ 
for a technical reason. 
If we pull back $e_1$ 
by a holonomy homomorphism 
$f:\pi_1S^1\to \Mg{}_{,*}$ of a surface bundle 
over the circle, 
it is clearly vanishing because $H^2(S^1,\Z)=0$. 
However 
Kitano showed in \cite{[Ki2]} that $e_1/48$, 
which depends on the spin structure of $\Sg$, 
makes sense as a bounded cohomology class 
in $H_b^2(S^1,\Z)\cong\R/\Z$ 
and it is given by the Rochlin invariant $\mu$ 
if the image of $f$ is contained in 
the Torelli group 
$\Ig{}_{,*}$. 
Combining Theorem 4.1 with the result of Kitano, 
we have the following on the level $2$ subgroup 
$\Mg{}_{,*}(2)
=\mathrm{Ker}\{\Mg{}_{,*}\to \mathrm{Sp}(2g,\Z/2)\}
\supset \Ig{}_{,*}$. 
That is, 
for $f:\Z \to \Mg{}_{,*}(2)$, 
the pull-back $f^*e_1/48\in H_b^2(\Z,\Z)$ 
is represented by 
$\mu(M_f,\tilde{\alpha})-\rho^{(2)}(\hat{M}_f)/16\in H^1(\Z,\R/\Z)$. 
The correspondence between them is given by 
the isomorphism 
$H_b^2(\Z,\Z)\cong H_b^1(\Z,\R/\Z)\cong 
H^1(\Z,\R/\Z)$. 

Here 
let us review the definition of the Rochlin invariant briefly. 
Let $(M,\alpha)$ be an oriented spin $3$-manifold with 
a spin structure $\alpha$. 
Then 
there exists a compact oriented spin $4$-manifold 
$(W,\beta)$ such that 
$\partial W=M$ and $\beta|_M=\alpha$. 
The Rochlin invariant 
$\mu(M,\alpha)\in \Q/\Z$ is defined by 
$$
\mu(M,\alpha)
=\frac{\sigma(W)}{16}
\quad 
\mathrm{mod}~ \Z
$$
and 
it does not depend on the choice of 
the $4$-manifold $W$ 
by virtue of Rochlin's theorem. 

Now 
Miller and Lee show in \cite{[ML]} that 
the Rochlin invariant of a spin $3$-manifold is 
a spectral invariant. 
To be more precice, 
it is described by the $\eta$-invariants as follows. 
Let $W,M$ be as in Section 3 and assume further that 
$W$ (hence also $M$) is a spin manifold. 
Let 
$\mathcal{D}$ denote the Dirac operator of $M$ 
acting on the spinor fields. 
This is a self-adjoint elliptic operator, 
so that 
the $\eta$-invariant $\eta_{\mathcal{D}}(M)$ 
is defined. 
We then have 
$$
\mathrm{ind}(\mathcal{D})
=
-\frac{1}{24}\int_Wp_1-\frac12\left\{\hbar+\eta_{\mathcal{D}}(M)\right\},
$$
where 
$\hbar$ is the dimension of the space of harmonic spinors on $M$. 
Combining the above formula and the index theorem mentioned 
in Section 3, 
we get 
$$
\sigma(W)+8~\mathrm{ind}(\mathcal{D})
=
-\eta(M)-4\left\{\hbar+\eta_{\mathcal{D}}(M)\right\}.
$$
The basic spin representations 
$\mathcal{S}^\pm$ of $\mathrm{Spin}(4)$ are quaternionic 
and 
hence the index of the Dirac operator is always even. 
We therefore obtain 
$$
\mu(M,\alpha)
=
-\frac{1}{16}\eta(M)-\frac14\left\{\hbar+\eta_{\mathcal{D}}(M)\right\}
\quad 
\mathrm{mod}~\Z,
$$
and 
it shows that 
the Rochlin invariant is a spectral invariant. 

Let us fix a spin structure $\alpha$ of $\Sg~(g\geq2)$. 
For each automorphism $f\in\Mg{}_{,*}(2)$, 
there exists the uniquely defined 
spin structure $\tilde{\alpha}$ on $M_f$ 
such that 
the restriction to each fiber is $\alpha$ and 
to the $S^1$-orbit of $*\in\Sg$ is 
the bounding spin structure. 
Applying the above formula 
to $(M_f,\tilde{\alpha})$ and 
substituting it for our description of 
the first Morita-Mumford class, 
we can conclude that 
$$
f^*e_1/48
=
-\frac{1}{16}\eta^{(2)}(\hat{M}_f)
-\frac14\left\{\hbar+\eta_{\mathcal{D}}(M_f)\right\}
\quad
\mathrm{mod}~\Z
$$
holds 
for an automorphism $f\in\Mg{}_{,*}(2)$. 

\vspace{2mm}

\noindent
\textit{Acknowledgements.} 
This research is partially supported by the Grant-in-Aid for Scientific Research
(No. 17740032), the Ministry of Education, Culture, Sports, Science and Technology, Japan.


\end{document}